\documentclass{amsart}

\usepackage[english]{babel}

\usepackage[a4paper,top=2cm,bottom=2cm,left=3cm,right=3cm,marginparwidth=1.75cm]{geometry}
\usepackage{mathrsfs}
\usepackage{amsmath, amssymb}
\usepackage{graphicx}
\usepackage[colorlinks=true, allcolors=blue]{hyperref}
\usepackage{mathtools}
\usepackage{amsfonts}
\usepackage{amsthm}
\usepackage{tikz-cd}

\newtheorem{thm}{Theorem}[section]
\newtheorem{intro}{Theorem}

\newtheorem{lem}[thm]{Lemma}
\newtheorem{prop}[thm]{Proposition}
\newtheorem{defn}[thm]{Definition}
\newtheorem{rem}[thm]{Remark}

\numberwithin{equation}{section}
\newcommand{\mycomment}[1]{}

\newcommand{\cL}{\mathcal{L}}

\newcommand{\tC}{\Tilde{C}}

\newcommand{\cO}{\mathcal{O}}

\newcommand{\unnumberedfootnote}[1]{%
  \begingroup
  \renewcommand{\thefootnote}{}%
  \footnotetext[0]{#1}%
  \endgroup
}

\usepackage[textwidth=0.8in]{todonotes}
\DeclareMathOperator{\Aut}{Aut}
\DeclareMathOperator{\id}{id}

\DeclareMathOperator{\Sym}{Sym}
\DeclareMathOperator{\Ker}{Ker}
\DeclareMathOperator{\im}{Im}
\DeclareMathOperator{\Pic}{Pic}
\DeclareMathOperator{\Nm}{Nm}
\DeclareMathOperator{\Ext}{Ext}
\DeclareMathOperator{\Cliff}{Cliff}
\DeclareMathOperator{\gon}{gon}
\DeclareMathOperator{\ev}{ev}
\DeclareMathOperator{\rk}{rank}

\title{Generic Prym-Torelli theorem for cyclic coverings}

\author[P. Bor\'owka]{Pawe\l{} Bor\'owka}
\address{Pawe\l{} Bor\'owka \newline Institute of Mathematics, Jagiellonian University in Krak\'ow\\ ul. prof. Stanisława Łojasiewicza 6, 
30-348 Kraków, Poland}
\email{pawel.borowka@uj.edu.pl}

\author[A. Shatsila]{Anatoli Shatsila}
\address{Anatoli Shatsila \newline Doctoral School of Exact and Natural Sciences, Jagiellonian University in Krak\'ow \\ ul. prof. Stanisława Łojasiewicza 6, 
30-348 Kraków, Poland}
\email{anatoli.shatsila@doctoral.uj.edu.pl}

\begin{document}
\maketitle

\unnumberedfootnote{\textit{2020 Mathematics Subject Classification:} 14H40, 14H30.}
\unnumberedfootnote{\textit{Key words and phrases:} Prym variety, Prym map, coverings of curves, Prym-Torelli, generic injectivity}

\begin{abstract}
    We prove the generic injectivity of the Prym map for étale cyclic coverings of degree \(n\geq 3\) of smooth curves of genus \(g\geq 7\). The proof is based on a description of the rank-one locus of the dual of the multiplication map $H^0(K_C\otimes\eta)\otimes H^0(K_C\otimes\eta^{-1}) \to H^0(K_C^{\otimes2})$ associated with a torsion line bundle $\eta$ defining the cover of $C$. In particular, we show that this locus recovers the bicanonical model of the base curve $C$ together with the linear system $|K_C\otimes \eta|$.
\end{abstract}

\section{Introduction}

The classical Torelli theorem states that a smooth curve is determined by its principally polarised Jacobian. For a Prym variety, the analogous question is whether one can recover the covering from the induced polarised abelian variety. Such results is called Prym-Torelli theorem. In the case of \'etale double covers, the Prym map is known to be generically injective for curves of genus at least $7$, although it is not injective on the whole moduli space (see \cite{BOprym} and the references therein). It is natural to ask the same question for cyclic coverings of higher degrees. It is worth noting that for $g\geq 3$, apart from some low degree cases with negative answers \cite{Faber, Xiao, ShatsilaTheta}, there has been no result showing Prym-Torelli theorem. For $g=2$, and more generally for hyperelliptic curves, only very recently, generic and global Prym-Torelli theorems have been shown for $d \geq 6$, see \cite{SG, NOPS, BNOS, ShatsilaTheta}.

Before going into the main results, let us introduce some notation. Let $f:\tC\to C$ be a connected \'etale cyclic cover of degree $n\geq 3$ of a smooth curve $C$ of genus $g\geq 2$. Such a cover is determined by a cyclic subgroup $\langle\eta\rangle\subset JC[n]$ of order $n$. By the Riemann-Hurwitz formula, $g(\tC)=n(g-1)+1$. Its Prym variety is defined as the connected component of the kernel of the norm map,
\[
P(f):=\Ker(\Nm_f:J\tC\to JC)^0.
\]
The restriction of the canonical principal polarisation of $J\tC$ induces a polarisation $\Xi_f$ on $P(f)$ of type $(1,\ldots,1,n,\ldots,n)$, where $n$ appears $g-1$ times. Thus, we obtain the Prym map
\[
\mathcal P_g[n]:\mathcal R_g[n]\longrightarrow
\mathcal A_{(n-1)(g-1)}^{(1,\ldots,1,n,\ldots,n)},
\qquad [f]\longmapsto [(P(f),\Xi_f)],
\]
where $\mathcal R_g[n]$ parametrizes connected \'etale cyclic covers of degree $n$ of smooth genus $g$ curves.


In \cite{BNOS}, it was conjectured that $\mathcal P_g[n]$ is generically injective for all $g\geq 3$ and $n\geq 4$. Our paper establishes generic injectivity of $\mathcal P_g[n]$ for all $g\geq 7$ and arbitrary $n\geq 3$, thus proving the conjecture in this range.

\begin{intro}[Theorem \ref{thm: main}]
For every $g\geq 7$ and $n\geq 3$, the Prym map $\mathcal P_g[n]$ is generically injective, i.e., generic Prym-Torelli theorem holds.
\end{intro}

The proof is based on the careful investigation of the codifferential of the Prym map. We first adapt the reconstruction result of \cite[Proposition~2.1]{BNOS} to cyclic coverings of arbitrary curves. More precisely, we show that for a non-hyperelliptic base curve, the automorphisms of the Prym variety induced from the deck group of the cover can be recovered from the polarised Prym itself as the subgroup of polarised automorphisms acting trivially on the kernel of the polarisation. This allows us to distinguish the character eigenspaces of $H^0(P(f),\Omega^1_{P(f)})$ intrinsically.

The main ingredient of the proof is a description of the rank-one locus of the dual of the multiplication map
\[
H^0(K_C\otimes\eta)\otimes H^0(K_C\otimes\eta^{-1})
\longrightarrow H^0(K_C^{\otimes2}).
\]
We prove that for a general curve of genus $g\geq7$ this locus consists precisely of the Schiffer variations of $C$. In particular, it gives the bicanonical model of the curve. To show this, we generalize the ideas of \cite{LangeOrtega} and \cite{MP}. The two factors of the rank-one tensors then recover the linear systems $|K_C\otimes\eta|$ and $|K_C\otimes\eta^{-1}|$, and hence the covering itself. 

The structure of the paper is as follows. We recall the relevant facts about automorphisms of Prym varieties and the codifferential of the Prym map in Section~2. Section~3 is devoted to the multiplication map and its rank-one locus. In Section~4 we prove the generic Prym-Torelli theorem.

\subsection*{Acknowledgements}

Anatoli Shatsila has been supported by the Polish National Science Centre project no. 2024/53/N/ST1/01634.
Pawe\l{} Bor\'owka has been supported by the Polish National Science Centre project number 2024/54/E/ST1/00330.

\section{Preliminaries}

Let $[f: \tilde{C} \to C] \in \mathcal{R}_g[n]$ and $\sigma\in\Aut(\tC)$ be a generator of the deck group, and denote by $(P, \Xi) = (P(f), \Xi_f)$ the Prym variety of $f$.  We will also think of $P$ as the complementary abelian subvariety to $f^*JC$ in $J\tC$, see \cite{BL}.
Since $\sigma$ preserves both $f^*JC$ and $P$, we abuse the notation and denote by $\sigma$ the restriction $\sigma|_P$ being a polarised automorphism of $(P, \Xi)$. We write $K(\Xi)=\Ker(\lambda_\Xi\colon P\to P^\vee)$ for the kernel of the polarising isogeny $\lambda_\Xi$, see \cite{BL}. The first observation on the way to showing that $\mathcal{P}_g[n]$ is generically injective is the following proposition, which is a counterpart of \cite[Proposition 2.1]{BNOS}. 

\begin{prop}
\label{prop: group}
    For the Prym variety $(P,\Xi)$ of a covering $f:\tC\to C$, define \begin{equation}
        G_P := \{\varphi \in \Aut(P, \Xi) \ : \ \varphi|_{K(\Xi)} = \id|_{K(\Xi)}\}.
    \end{equation}
    If $C$ is non-hyperelliptic then
    $
    G_P = \langle \sigma \rangle \simeq \mathbb{Z}/n\mathbb{Z}.
    $
    In particular, the group $\langle \sigma \rangle$ is intrinsic in $\Aut(P, \Xi)$.
\end{prop}
\begin{rem}
    Note that \cite[Proposition 2.1]{BNOS} treats the case when $C$ is hyperelliptic and then $G_P= \langle \sigma, -j\rangle \simeq D_{n}$ where \(D_n\) is the dihedral group with $2n$ elements and $j|_P$ is an automorphism of $(P, \Xi)$ induced by a lift $j \in \Aut(\tilde{C})$ of the hyperelliptic involution $\iota$ on $C$. This shows in particular that $C$ being non-hyperelliptic is a condition satisfied for all elements (or none) in the fibre of the Prym map. 
\end{rem}
\begin{proof}[Proof of Proposition 2.1]
    Let $f:\tC\to C$ be a given cover. Since $f^*JC$ and $P$ are complementary with respect to the principal polarisation on $J\tC$, the addition map
\[
 m:f^*JC\times P\longrightarrow J\tC
\]
is an isogeny with kernel $\{(x,-x):x\in f^*JC\cap P\}$, where $f^*JC\cap P=K(\Xi)$, see \cite{BL}. Let $\varphi\in G_P$. Since $\varphi$ is the identity on $K(\Xi)$, the map $\id_{f^*JC}\times\varphi$ descends to an automorphism $\widetilde\varphi$ of $J\tC$, which acts as the identity on $f^*JC$ and as $\varphi$ on $P$.

We claim that $\widetilde\varphi$ preserves the canonical principal polarisation on $J\tC$. Indeed, the pullback $m^*[\Theta_{\tC}]$ is the product of the restricted polarisations on $f^*JC$ and $P$, hence is preserved by $\id_{f^*JC}\times\varphi$. Since the pullback by an isogeny is injective on N\'eron--Severi groups, we obtain $\widetilde\varphi^*[\Theta_{\tC}]=[\Theta_{\tC}]$. Therefore, by the strong Torelli theorem, we can write
\[
 \widetilde\varphi=\varepsilon a_*,\qquad a\in\Aut(\tC),\quad \varepsilon\in\{\pm1\}.
\]

We claim that $a$ descends to an automorphism of $C$. Since $\widetilde\varphi$ is the identity on $f^*JC$ and preserves $P$, we have $\Nm_f\circ(\widetilde\varphi-\id)=0$. Thus, for any $p,q\in\tC$ with $f(p)=f(q)$, we get
\[
 [f(a(p))-f(a(q))]=0\quad\text{in }JC.
\]
The Abel-Jacobi map of $C$ is injective, so $f(a(p))=f(a(q))$. It follows that $a$ descends to an automorphism $a_0$ of $C$. Moreover,
\[
 f^*\circ(a_0)_*=a_*\circ f^*
=\varepsilon\widetilde\varphi\circ f^*=\varepsilon f^*.
\]
Since $f^*$ has finite kernel, we conclude that $(a_0)_*=\varepsilon\id_{JC}$.

If $\varepsilon=1$, then $a_0=\id_C$, so $a=\sigma^k$ for some $k$. If $\varepsilon=-1$, then $C$ is hyperelliptic that is a contradiction. Therefore, every element of $G_P$ is of the form $\sigma^k$.

Conversely, $\sigma$ preserves $\Xi$ and acts trivially on $f^*JC$, hence also on $K(\Xi)=f^*JC\cap P$, so it belongs to $G_P$. Moreover, its action on $P$ is faithful: $\sigma$ acts diagonally on the non-trivial character eigenspaces of $H^0(P,\Omega_P^1)$.
\end{proof}
\begin{rem}\label{singP}
    Proposition \ref{prop: group} makes a difference between $n=2$ and $n\geq 3$ cases because in $n=2$ case, the deck automorphism $\sigma$ restricts to $-1$ on $P$ and gives no additional information. For $n\geq 3$, we know that the Prym variety is an abelian variety with an automorphism of order $n$.
\end{rem}

\subsection{The codifferential of the Prym map}

In this subsection, we give a description of the codifferential of the map $\mathcal{P}_g[n]$, following \cite{LangeOrtega}. Throughout the paper we use the multiplicative convention for line bundles with $\eta^k=\eta^{\otimes k}$ and for a divisor D, we write $\eta(D)=\eta\otimes \cO(D)$. We will also use the notions and identifications of tangent and cotagent spaces. Unfortunately, as pointed out in Remark \ref{singP}, the Prym varieties contain additional automorphisms, so they lie in the singular locus of the coarse moduli space of polarised abelian varieties. One can treat this issue in 3 different ways. Firstly, one can use a modern language of stacks or equivalently, orbifolds that are smooth. Secondly, more classically, one can put a level structure to make the moduli space fine and consider lifted maps. Lastly, analytically, one can lift a problem to a sufficiently small disk in the Siegel space. We leave to the reader any technical details that come from their preferred point of view.


Let $[f\colon\tilde C\to C]\in\mathcal R_g[n]$, let $\sigma$ be a generator of the deck group and $\zeta$ a primitive $n$-th root of unity. The group $\langle\sigma\rangle$ acts on $f_*\mathcal O_{\tilde C}$, which splits as
$$f_*\mathcal O_{\tilde C}=\bigoplus_{i=0}^{n-1}\cL_i,\qquad \sigma^*|_{\cL_i}=\zeta^i\id,$$
where each $\cL_i$ is a line bundle, $\cL_0=\mathcal O_C$ and $\cL_i\cong \cL_1^{i}$. We put $\eta:=\cL_1$ and therefore we fix a generator of the subgroup of $JC[n]$ defining $f$. Changing $\sigma$ or $\zeta$ replaces $\eta$ by another generator. Using the fact that f is \'etale and the projection formula, we obtain
\begin{equation}
K_{\tC} = f^*K_C, \qquad H^0(K_{\tilde C})=\bigoplus_{i=0}^{n-1}V_i,\ \text{ where }\ V_i:=H^0(K_C\otimes\eta^i)=\Ker(\sigma^*-\zeta^i\id).
\end{equation}

The summand $V_0=f^*H^0(K_C)$ can be identified with the cotangent space of $f^*JC$, hence
$$(T_0P)^*\cong H^0(\Omega^1_P)\cong\bigoplus_{i=1}^{n-1}V_i .$$

In particular, this decomposition is determined by the group $\langle\sigma\rangle\subset\operatorname{Aut}(P,\Xi)$, up to the choice of the generator $\sigma$ and of $\zeta$.


The codifferential of the Prym map at $[f]$ is the projection of multiplication map, see \cite[Proposition~4.1]{LangeOrtega} $$(d_{(C,\eta)}\mathcal{P}_g[n])^*: \left(T_P\mathcal{A}_{(n-1)(g-1)}^{(1,\ldots,1,n,\ldots,n)}\right)^* \simeq \Sym^2H^0(\Omega_P^1) \to H^0(C, K_C^{2}) \simeq \left(T_{(C,\eta)}\mathcal{R}_g[n]\right)^{*}$$ 
that can be rewritten as
\begin{equation}\label{eq:prym-codiff}
    (d_{(C,\eta)}\mathcal{P}_g[n])^*: \bigoplus_{i=1}^{n-1} \Sym^2 V_i \oplus \bigoplus_{i < j} (V_i \otimes V_j) \to H^0(C, K_C^{2}).
\end{equation}
Note that, only the $\langle \sigma \rangle$-invariant summands contribute to the map: summands $V_i\otimes V_j$ map non-trivially to $H^0(C,K_C^{2})$ only if $i+j\equiv 0\pmod n$ including $\Sym^2V_{\frac n2}$ if $n$ is even.


\section{Schiffer variations under twisted multiplication}

Let $C$ be a smooth projective curve of genus $g \geq 2$ and $p \in C$ be a point. Using the long exact sequence associated to 
$$0 \to T_C \to T_C(p) \to T_C(p)|_p \to 0$$ and the fact that $h^0(T_C(p)) = 0, h^0(T_C(p)|_p) = 1$ the image of the connecting map $$\delta_p: H^0(T_C(p)|_p) \to H^1(T_C)$$ is a line. We define a Schiffer variation as follows

\begin{defn}
    A non-zero element $\xi_p \in \im \delta_p$ is called a Schiffer variation at $p$. 
\end{defn}

We use a classical fact that for $g \geq 3$ the map $$C \ni p \mapsto [\xi_p] \in \mathbb PH^1(T_C) \simeq \mathbb{P}H^0(K_C^{2})^*$$ is the bicanonical embedding of $C$. 

Now, let $\eta \in \Pic^0(C)$ be a non-trivial line bundle on $C$, and denote by $V = V_1= H^0(K_C \otimes \eta)$ and $W = V_{n-1}= H^0(K_C\otimes \eta^{-1})$. Consider the multiplication map $$\mu_{\eta}: V \otimes W \to H^0(K_C^{2})$$ and its dual $$\mu_{\eta}^*: H^0(K_C^{2})^* \to V^*\otimes W^*.$$
Not that the map $\mu_{\eta}$ is the restriction of the Prym codifferential to the summand $V\otimes W$ in \eqref{eq:prym-codiff}.

The main goal of this section is to relate the rank-one locus inside $V^*\otimes W^*$ with the image of $\mu^*_{\eta}$. We start with a few auxiliary results that we need in the proof of Theorem \ref{thm: rank1locusisnice}.

\begin{lem}
\label{lem: torsion-avoidance}
    Let $C$ be a general smooth genus $g \geq 2$ curve and $k := \lfloor \frac{g-1}{2}\rfloor$.  Then $$C_k - C_k = \{\mathcal{O}_C(E - D): \  E, D \in C^{(k)}\}\subset \Pic^0(C)$$ contains no non-zero torsion points.
\end{lem}

\begin{proof}
    Suppose that $\mathcal{O}_C(D-E) \in C_d-C_d$ is a non-trivial torsion line bundle, where $D$ and $E$
    are effective divisors of degree $d$. Cancelling their common part, we can assume that $D$ and $E$ have disjoint supports and degree $e \leq d$. It follows from \cite{Kempf} that
    $h^0(K_C(-D-E))=0$. By Riemann-Roch,
    \[
        0=h^0(K_C(-D-E))
        =g-1-2e+h^0(D+E)>g-1-2e,
    \]
    so $d \geq e > \frac{g-1}{2}$, and the result follows.
\end{proof}

\begin{lem}
\label{lem: twisted-clifford}
    Let $C$ be a general smooth genus $g \geq 7$ curve, and let $\eta \in \Pic^0(C)$ be a non-trivial line bundle of order $n\geq 3$. Then for any effective divisor $D$ of degree $d$ on $C$ satisfying $2 \leq d \leq 2g-2-\gon(C)$ we have $$h^0(\eta(D)) + h^0(\eta^{-1}(D)) \leq d - 2.$$
\end{lem}

\begin{proof}
    Firstly, note that $h^0(\eta(D)) = h^0(\eta^{-1}(D)) = 0$ for any effective divisor $D$ of degree at most 3, which establishes the result for $d=2, 3$. This is because if $h^0(\eta(D)) \geq 1$ then $\eta=\cO(E-D)$ that contradicts Lemma \ref{lem: torsion-avoidance}.
    
    From now on we assume that $d \geq 4$. We also assume $g \geq 9$. Under this assumption, we have $\Cliff(C) \geq 4$. We will prove that $2h^0(\eta(D)) \leq d-2$ and $2h^0(\eta^{-1}(D)) \leq d-2$, which implies the assertion. Let us prove this for the line bundle $\eta(D)$; the proof for $\eta^{-1}(D)$ is analogous. 
    
    \emph{Case 1.} If $h^0(\eta(D)) \leq 1$, then the assertion is automatic. 

    \emph{Case 2.} If $h^0(\eta(D)) \geq 2$ and $h^1(\eta(D)) \geq 2$ then $$d - 2h^0(\eta(D)) + 2 \geq \Cliff(C) \Longrightarrow 2h^0(\eta(D)) \leq d+2-\Cliff(C) \leq d-2.$$

    \emph{Case 3.} If $h^1(\eta(D)) \leq 1$, then Riemann-Roch implies $$2h^0(\eta(D)) = 2d - 2g + 2 + 2h^1(\eta(D)) \leq 2d-2g+4 \leq d - \Cliff(C) \leq d-2.$$
    Therefore, the result is true for any $g \geq 9$. 
    
    Assume that $g = 7$ or $g = 8$, so that $\Cliff(C) = 3$ and $\gon(C) = 5$, in particular $4 \leq d \leq 2g-7$. The same reasoning as above shows that $2h^0(\eta(D)) \leq d-1$ and $2h^0(\eta^{-1}(D)) \leq d-1$, so we only need to rule out the case $2h^0(\eta(D)) = 2h^0(\eta^{-1}(D)) = d-1$. Since $d$ is odd, we have $(g,d) \in \{(7,5), (7,7), (8,5), (8,7), (8,9)\}$. The case $(g,d) = (8,7)$ is impossible since there is no $g^2_7$ on a general genus 8 curve.

    Let us show that in all remaining cases $C$ possesses two distinct $g^1_5$'s $L_1$ and $L_2$ such that $L_1 \otimes L_2^{-1} \simeq \eta^{2}.$ 
    If $d = 5$ then $L_1 = \eta(D)$ and $L_2 = \eta^{-1}(D)$ work. If $(g,d) \in \{(7,7), (8,9)\}$ then we can take $L_1 = K_C\otimes \eta(-D)$ and $L_2 = K_C\otimes \eta^{-1}(-D)$. It follows from Lemma \ref{lem: two_g^1_5} below that we get a contradiction with $\eta^2$ being torsion, so the proof is finished.
\end{proof}

\begin{lem}
\label{lem: two_g^1_5}
Let $C$ be a general curve of genus $7$ or $8$. Then $W^1_5-W^1_5$ contains no non-zero torsion points.
\end{lem}

\begin{proof}
Suppose $L_1\otimes L_2^{-1}$ is a non-trivial torsion line bundle, where $L_1,L_2\in W^1_5$. Since $C$ is general, $h^0(L_i)=2$, and since $\mathrm{gon}(C)=5$, both $L_i$ are base-point free. Let $f_i:C\to\mathbb P^1$ be the morphisms given by $|L_i|$, and set $$f=(f_1,f_2):C\to\mathbb P^1\times\mathbb P^1.$$

The morphism $f$ is birational onto its image. Otherwise its degree onto the image divides $\deg f_1=5$, so it equals 5. Then $f(C)$ is a curve of bidegree $(1,1)$, whose projections to both factors are isomorphisms, and hence $L_1\cong L_2$, a contradiction.

Hence $f(C)$ has bidegree $(5,5)$, so $p_a(f(C))=16>g(C)$ and hence $f$ is not a closed immersion, so it fails to separate either two distinct points $p,q$ or a tangent vector at some point $p$. Accordingly, let $p+q$ or $2p$ (if $p=q$) be an effective divisor of degree 2 contained in a fibre of each $f_i$.
Since the fibres of $f_i$ are the divisors of $|L_i|$, we get $h^0(L_i(-p-q))\ge1$, that is, $L_i\cong\mathcal O_C(p+q+E_i)$ for some effective $E_i$ of degree 3. Therefore

$$L_1\otimes L_2^{-1}\cong\mathcal O_C(E_1-E_2)\in C_3-C_3,$$
which contradicts Lemma \ref{lem: torsion-avoidance}.

\end{proof}

Recall that $V = H^0(K_C \otimes \eta)$, $W = H^0(K_C\otimes \eta^{-1})$, the multiplication map $\mu_{\eta}: V \otimes W \to H^0(K_C^{2})$ and its dual is $$\mu_{\eta}^*: H^0(K_C^{2})^* \to V^*\otimes W^*.$$ In the proof of the following theorem, we follow the ideas of of the proof of \cite[Thm 2.8]{MP}.

\begin{thm}
\label{thm: rank1locusisnice}
    Let $C$ be a general smooth genus $g \geq 7$ curve, and let $\eta \in \Pic^0(C)[n]$ be a non-trivial line bundle of order $n \geq 3$. Then the dual to the multiplication map $\mu_{\eta}^*$ is injective and for every non-zero $\xi \in H^1(T_C)$ we have $$\mu_{\eta}^*(\xi) \text{ has rank } 1 \Leftrightarrow [\xi] = [\xi_p] \text{ for some }p\in C.$$
\end{thm}

\begin{proof}
    Let $0 \neq \xi \in H^1(T_C)$. Since $H^1(T_C) \simeq H^1(K_C^{-1}) \simeq \Ext^1(K_C\otimes \eta, \eta)$, $\xi$ determines a non-split extension of vector bundles
    \begin{equation}
    \label{extension}
        0 \to \eta \to F_{\xi} \to K_C\otimes \eta \to 0.
    \end{equation}
    Since $h^0(\eta) = 0$ and $h^0(K_C\otimes \eta) = g-1$, the long exact sequence in cohomology reads $$0 \to H^0(F_{\xi}) \to H^0(K_C\otimes \eta) \xrightarrow[]{\partial_{\xi}} H^1(\eta) \to \ldots$$ Moreover, since $W^* \simeq H^1(\eta)$, the connecting homomorphism $\partial_{\xi}$ (seen as an element of $V^*\otimes W^*$) is precisely $\mu^*_{\eta}(\xi)$, hence $$h^0(F_{\xi}) = g-1 - \rk \mu^*_{\eta}(\xi).$$  
    
    Assume that $\rk \mu^*(\xi) \leq 1$. Since $\det F_{\xi} = K_C\otimes \eta^2$ and $\eta^2 \neq \mathcal{O}_C$, we have $h^0(\det F_{\xi}) = g-1$. Since $g \geq 7$, $$2h^0(F_{\xi}) - 3 = 2g-5-2\rk \mu^*_{\eta}(\xi)>g-1 = h^0(\det F_{\xi}).$$
     Therefore, \cite[Lemma 3.9]{ParanjapeRamanan} implies that $F_{\xi}$ has a saturated line subbundle $A \subset F_{\xi}$ with $h^0(A) \geq 2$. The composition $$A \hookrightarrow F_{\xi} \to K_C\otimes \eta$$ must be non-zero, since otherwise $A \hookrightarrow \eta$, contradicting $h^0(\eta) = 0$. Therefore, $$A \simeq K_C\otimes \eta(-D)$$ for some effective divisor $D \neq 0$ (recall that the extension does not split). Since $h^0(A) \geq 2$, we have $\deg A \geq \gon(C)$ so $$1 \leq \deg D \leq 2g-2-\gon(C).$$ Taking determinants in the short exact sequence $$0 \to A \to F_{\xi} \to F_{\xi}/A \to 0$$ yields $F_{\xi}/A \simeq \eta(D)$, so we have an exact sequence $$0 \to K_C\otimes \eta(-D) \to F_{\xi} \to \eta(D) \to 0.$$ Taking global sections and applying Riemann-Roch to $K_C\otimes \eta(-D)$ yields $$\deg D - \rk \mu^*_{\eta}(\xi) \leq h^0(\eta(D)) + h^0(\eta^{-1}(D)).$$

    If $\rk \mu^*_{\eta}(\xi) = 0$ then Lemma \ref{lem: twisted-clifford} implies that $\deg D = 1$. Then, $h^0(\eta(D)) \geq 1$ yields $\eta \simeq \mathcal{O}_C(E-D)$ for effective divisors $E, D$ of degree 1, contradicting Lemma \ref{lem: torsion-avoidance}, and similar for $h^0(\eta^{-1}(D)) \geq 1$, so we reach a contradiction in the case $\rk \mu^*_{\eta}(\xi) = 0$. 
    Therefore, we have shown that $\mu^*_{\eta}$ is injective.

    If $\rk \mu^*_{\eta}(\xi) = 1$ then again Lemma \ref{lem: twisted-clifford} implies that $\deg D = 1$, i.e. $D = p$ and $A \simeq K_C\otimes \eta(-p)$. Thus, the pullback of Extension \ref{extension} along $K_C\otimes \eta(-p) \to K_C\otimes \eta$ splits, hence $$\xi \in \ker \left(H^1(T_C)\simeq \Ext^1(K_C\otimes \eta, \eta) \to \Ext^1(K_C\otimes \eta(-p), \eta) \simeq H^1(T_C(p))\right) \Longrightarrow [\xi] = [\xi_p].$$

    It follows from Lemma \ref{lem: torsion-avoidance} that $K_C\otimes \eta$ and $K_C\otimes \eta^{-1}$ are very ample, hence $\mu^*_{\eta}(\xi_p) = (\ev_p^V \otimes \ev_p^{W}) \in V^*\otimes W^*$ is a rank 1 tensor, which establishes another implication.
\end{proof}

\begin{rem}
    Injectivity of $\mu^*_{\eta}$ also follows from \cite[Theorem 1]{Butler} (see also \cite[Proposition 5.5]{LangeOrtega}).
\end{rem}

\begin{rem}
\label{rmk: theoremfailsforsmallgenera}
    The conclusion of Theorem \ref{thm: rank1locusisnice} fails for $g\leq 5$. Indeed, let $g=5$ and suppose that $\mu_\eta^*$ is injective. Since $h^0(K_C\otimes\eta)=h^0(K_C\otimes\eta^{-1})=4$ and $h^0(K_C^{\otimes 2})=12$, its projectivized image is a copy of $\mathbb P^{11}$ in $\mathbb P^{15}$. The rank-one tensors form the Segre variety $\mathbb P^3\times\mathbb P^3\subset\mathbb P^{15}$. Thus, by the projective dimension theorem,
    \[
        \dim\bigl(\mathbb P(\operatorname{im}\mu_\eta^*)\cap(\mathbb P^3\times\mathbb P^3)\bigr)
        \geq 11+6-15=2.
    \]
    In particular, this locus cannot be the bicanonical image of $C$. If $\mu_\eta^*$ is not injective, the theorem already fails. The same argument applies in genus $4$, while in genus $3$ the map $\mu_\eta^*$ cannot be injective for dimension reasons. We do not know whether Theorem \ref{thm: rank1locusisnice} fails for $g=6$.
\end{rem}

\section{Generic injectivity of the cyclic Prym map for \texorpdfstring{$g \geq 7$}{}}

The main goal of this section is to show generic injectivity of the Prym map $\mathcal{P}_g[n]$ for $g \geq 7$. 

\begin{thm}
\label{thm: main}
    Let $g \geq 7, n \geq 3$ be integers. The Prym map of cyclic covers of degree $n$ $$\mathcal{P}_g[n]: \mathcal{R}_g[n] \to \mathcal{A}_{(n-1)(g-1)}^{(1,\ldots,1,n,\ldots,n)}$$ is generically injective.
\end{thm}

\begin{proof}

    Let $Y$ be the closure of the image of $\mathcal P_g[n]$. By Theorem \ref{thm: rank1locusisnice}, for a general cover $(C,\langle\eta\rangle) \in \mathcal{R}_g[n]$ the multiplication map $$\mu_\eta:\ H^0(K_C\otimes\eta)\otimes H^0(K_C\otimes\eta^{-1})\longrightarrow H^0(K_C^{\otimes2})$$ is surjective. Since this map is a restriction of the codifferential of $\mathcal P_g[n]$, the Prym map is generically finite onto $Y$. In particular, $\dim Y=3g-3$.
    
    Let $y=[(P,\Xi)]$ be a general point of $Y$, and suppose that $$(C_1,\langle\eta_1\rangle),\qquad (C_2,\langle\eta_2\rangle)$$ are two points in its fiber, which we can take to be general. We can assume that both curves are non-hyperelliptic, that $Y$ is smooth at $y$, and that the Prym map is unramified at these two points. We will show that the two covers are isomorphic.

    By Proposition~\ref{prop: group}, the group $G_P=\langle\sigma\rangle$ is determined by $(P,\Xi)$. Fix a generator $\sigma$ and a primitive $n$-th root of unity $\zeta$. Let $$V=\Ker(\sigma^*-\zeta\,\id),\qquad W=\Ker(\sigma^*-\zeta^{-1}\,\id)$$ be the corresponding eigenspaces in $H^0(P,\Omega_P^1)$. Choosing the generators of the two deck groups of the coverings to be compatible with $\sigma$, we have
    \begin{equation}\label{eq:prym-intrinsic-VW}
    V\simeq H^0(K_{C_j}\otimes\eta_j),\qquad W\simeq H^0(K_{C_j}\otimes\eta_j^{-1}),\qquad j=1,2.
    \end{equation}
    Thus the two multiplication maps $\mu_{\eta_j}$ have the same source $V\otimes W$. On the other hand, since the Prym map is unramified at these points and $\dim Y=3g-3$, its codifferentials induce isomorphisms
    \begin{equation}\label{eq:prym-two-cotangent-isos}
    \theta_j:T_y^*Y\xrightarrow{\ \sim\ } H^0(K_{C_j}^{\otimes2}),\qquad j=1,2.
\end{equation}
    Consider the cotangent map of the inclusion $Y \hookrightarrow \mathcal{A}_{(n-1)(g-1)}^{(1,\ldots,1,n,\ldots,n)}$, restricted to the summand $V\otimes W\subset\Sym^2 H^0(P,\Omega_P^1)$. It gives a map $$\alpha:V\otimes W\longrightarrow T_y^*Y.$$
    By the description of the Prym codifferential in \eqref{eq:prym-codiff}, we obtain
    \begin{equation}\label{eq:prym-common-tensor}
    \mu_{\eta_j}=\theta_j\circ\alpha,\qquad j=1,2.
    \end{equation}
    In particular, the two multiplication maps are identified through the cotangent space of $Y$.

    Consider the locus with reduced structure $$Z=\bigl\{[v]\in\mathbb P(T_yY): \rk(\alpha^*(v))=1\bigr\}_{\mathrm{red}}.$$
    By \eqref{eq:prym-common-tensor}, the dual map $\alpha^*$ can be identified with $\mu_{\eta_j}^*$ via $\theta_j^*$. Now, Theorem \ref{thm: rank1locusisnice} shows that $Z$ is the bicanonical image of $C_j$ for $j=1,2$. In particular,
    \begin{equation}\label{eq:prym-recover-curve}
    C_1\simeq Z\simeq C_2.
    \end{equation}

    Finally, every rank-one tensor $\alpha^*(v)$ determines a point of the Segre variety $$\mathbb P(V^*)\times\mathbb P(W^*)\subset\mathbb P(V^*\otimes W^*),$$
    and hence two morphisms $$Z\longrightarrow\mathbb P(V^*),\qquad Z\longrightarrow\mathbb P(W^*).$$
    Under the identification $Z\simeq C_j$, the formula $$\mu_{\eta_j}^*(\xi_p)=\ev_p^{K_{C_j}\otimes\eta_j}\otimes \ev_p^{K_{C_j}\otimes\eta_j^{-1}}$$
    shows that these morphisms are defined by complete linear systems $|K_{C_j}\otimes\eta_j|$ and $|K_{C_j}\otimes\eta_j^{-1}|$, respectively. Thus, the pullback of $\mathcal O_{\mathbb P(V^*)}(1)$ by the first morphism is $K_{C_j}\otimes\eta_j$. Since the canonical bundle of $Z$ is already known, we recover $\eta_j$ as well. Consequently the two cyclic subgroups $\langle\eta_j\rangle$ agree, so the covers are isomorphic.
\end{proof}

\begin{rem}
    By \cite[Lemma 9.1]{Faber}, for $n=3$ the Prym map is never injective, so generic injectivity cannot be strengthened to global injectivity. For $n\geq4$ we do not know any examples of non-injectivity locus. 
\end{rem}

\subsection*{Open questions}
It would be interesting to extend our results to curves of smaller genera. Our proof uses the assumption $g\geq 7$ in an essential way. The proof of Lemma~\ref{lem: twisted-clifford} uses the fact that a general curve has Clifford index at least $3$, together with the absence of non-trivial torsion points in $C_3-C_3$. Neither of these properties holds for a general curve of genus less than $7$.

For $g=3,4,5$, Remark \ref{rmk: theoremfailsforsmallgenera} shows that one piece of the multiplication map cannot recover the curve from its rank-one locus. One may instead try to use the full codifferential of the Prym map, taking into account all opposite-character summands, or to study higher-rank loci. 

Another interesting direction is to investigate whether the Prym maps are injective, not only generically. Our result, together with \cite{BNOS}, restricts potential non-injectivity locus, however the ideas of the paper cannot be used to prove global injectivity. On the other hand, we have not been able to generalize the trigonal construction of Faber \cite[\S 9]{Faber} to higher degrees.


\subsection*{AI disclosure}

The proof of the main theorem was developed with the assistance of ChatGPT Astra 6. The authors refined and simplified the proof of Theorem \ref{thm: rank1locusisnice} suggested by the model, verified all mathematical details, and take full responsibility for the correctness of the results.

\bibliographystyle{alpha}
\bibliography{sample}

\end{document}